\documentclass[fleqn]{mat01}
\usepackage{times,mathtimy,amssymb,latexsym}
\begin{document}

\def\d{\mbox{\rm d}}
\def\e{\mbox{\rm e}}


\def\mand{\quad \mbox{and} \quad}
\def\cS{{\mathcal S}}
\def\cR{{\mathcal R}}

\newtheorem{theore}{Theorem}
\renewcommand\thetheore{\arabic{section}.\arabic{theore}}
\newtheorem{theor}{\bf Theorem}
\newtheorem{lem}[theore]{\it Lemma}
\newtheorem{propo}[theore]{\rm PROPOSITION}
\newtheorem{coro}[theore]{\rm COROLLARY}
\newtheorem{definit}[theore]{\rm DEFINITION}
\newtheorem{probl}[theore]{\it Problem}
\newtheorem{exampl}[theore]{\it Example}

\title{Arithmetic properties of the Ramanujan function}

\markboth{Florian Luca and Igor E Shparlinski}{Arithmetic properties of the Ramanujan function}

\author{FLORIAN LUCA$^{1}$ and IGOR E SHPARLINSKI$^{2}$}

\address{$^{1}$Instituto de Matem{\'a}ticas, Universidad Nacional Aut\'onoma de
M{\'e}xico, C.P.~58089, Morelia, Michoac{\'a}n, M{\'e}xico\\
\noindent $^{2}$Department of Computing, Macquarie University, Sydney, NSW~2109, Australia\\
\noindent E-mail: fluca@matmor.unam.mx; igor@ics.mq.edu.au}

\volume{116}

\mon{February}

\parts{1}

\pubyear{2006}

\Date{MS received 2 December 2004\\[1.6pc]
\noindent {\it Dedicated to T~N~Shorey on his sixtieth
birthday}\vspace{.6pc}}

\begin{abstract}
We study some arithmetic properties of the Ramanujan function $\tau(n)$,
such as the largest prime divisor $P(\tau(n))$ and the number of
distinct prime divisors $\omega(\tau(n))$ of $\tau(n)$ for various
sequences of $n$. In particular, we show that \hbox{$P(\tau(n)) \geq (\log
n)^{33/31 + o(1)}$} for infinitely many $n$, and
\begin{equation*}
P(\tau(p)\tau(p^2)\tau(p^3)) > (1+o(1))\frac{\log\log p\log\log\log p}
{\log\log\log\log p}
\end{equation*}
for every prime $p$ with \hbox{$\tau(p)\neq 0$}.
\end{abstract}


\maketitle

\section{Introduction}

Let $\tau(n)$ denote the {\it Ramanujan function} defined by the
expansion
\begin{equation*}
X\prod_{n =1}^\infty (1-X^n)^{24} = \sum_{n=1}^\infty \tau(n) X^n,
\qquad |X| < 1.
\end{equation*}

For any integer $n$ we write $\omega(n)$ for the number of
distinct prime factors of $n$, $P(n)$ for the largest prime factor
of $n$ and $Q(n)$ for the largest square-free factor of $n$ with
the convention that \hbox{$\omega(0)=\omega(\pm 1) = 0$} and
\hbox{$P(0)=P(\pm 1) = Q(0)=Q(\pm 1)=1$}.

In this note, we study the numbers  $\omega(\tau(n))$, $P(\tau(n))$ and
$Q(\tau(n))$
as $n$ ranges over various sets of positive integers.

The following basic properties of $\tau(n)$
underline  our approach which is similar to those of~\cite{MMS,Shorey}:

\begin{itemize}
\item $\tau(n)$ is an integer-valued multiplicative function; that is,
\hbox{$\tau(m)\tau(n) = \tau(mn)$} if \hbox{$\gcd(m,n)=1$}.

\item For any prime $p$, and an integer \hbox{$r \geq 0$},
$\tau(p^{r+2}) = \tau(p^{r+1})\tau(p) - p^{11} \tau(p^{r})$,
where $\tau(1) = 1$.
\end{itemize}

In particular, the identity
\begin{equation}\label{eq:quadr}
\tau(p^2) = \tau(p)^2 - p^{11}
\end{equation}
plays a crucial role in our arguments.

It is also useful to recall that by the famous result of Deligne
\begin{equation}
\label{eq:Deligne}
|\tau(p)|\leq 2 p^{11/2} \mand |\tau(n) |  \leq n^{11/2 + o(1)}
\end{equation}
for any prime $p$ and positive integer $n$ (see~\cite{Murty}).

One of the possible approaches to studying arithmetic properties
of $\tau(n)$ is to remark that the values $u_r =\tau(2^r)$ form a Lucas
sequence satisfying the following binary recurrence relation
\begin{equation}
\label{eq:Bin Rec}
u_{r+2} = -24 u_{r+1}- 2048 u_r, \qquad r = 0,1,\dots,
\end{equation}
with the initial values  $u_0 = 1$, $u_1 = -24$.
By the primitive divisor theorem for Lucas sequences which claims that
each  sufficiently large  term $u_r$ has at least one new prime
divisor (see~\cite{BHV} for the most
general form of this assertion), we conclude that
\begin{equation*}
\omega\left(\prod_{r\leq z}\tau(2^r)\right)\geq z+O(1),
\end{equation*}
leading to the inequality
\begin{equation*}
\omega\left(\prod_{\substack{n\leq x\\ \tau(n)\neq 0}}
\tau(n)\right)\geq \left(\frac{1}{\log 2}+o(1)\right)\log x
\end{equation*}
as $x\rightarrow \infty$. In particular, we derive that for
infinitely many $n$,
\begin{equation*}
P(\tau(n)) \geq \log n \log \log n.
\end{equation*}
A stronger conditional result, under the $ABC$-conjecture, is given
in~\cite{MuWo}. We also have
\begin{equation*}
Q(\tau(n)) \geq  n^{(\log 2 + o(1))/ \log \log \log n}
\end{equation*}
for infinitely many $n$ (see eq.~(16) in~\cite{Stew}).

Furthermore, since $u_r | u_s $, whenever $r+1|s+1$, it follows
that if for sufficiently large $s$ we set $k={\rm lcm}[2, \dots,
s+1]-1$, then $\tau(2^k)$  is divisible by $\tau(2^r)$ for all
\hbox{$r\leq s$}. Thus, setting $n =2^k$ we get
\begin{align*}
\omega(\tau(n))\geq  s + O(1) = \left(\frac{1}{\log
2}+o(1)\right)\log k\ge \left(\frac{1}{\log 2}+o(1)\right)\log\log
n
\end{align*}
as $n\rightarrow \infty$. Here, we  use different approaches to
improve on these bounds.

Our results are based on some bounds for smooth numbers, that is,
integers $n$ with restricted $P(n)$ (see~\cite{HiTe,Te}). We also
use results on $\cS$-unit equations (see~\cite{Eve}). We recall
that for a given finite set of primes $\cS$, a rational $u
=s/t\neq 0$ with $\gcd(s,t) =1$ is called an $\cS$-unit if all
prime divisors of both $s$ and $t$ are contained in $\cS$.
Finally, we also use bounds on linear forms in $q$-adic
logarithms (see~\cite{Yu}).

We recall that in~\cite{MM} it is shown under the extended Riemann
hypothesis  that $\omega(\tau(p))\sim \log\log p$  holds for
almost all primes $p$ and that $\omega(\tau(N))\sim 0.5(\log\log
N)^2$ holds for almost all positive integers $N$.

Throughout the paper, the implied constants in the symbols `$O$',
`$\gg$' and `$\ll$' are absolute (recall that the notations $U \ll
V$ and $V\gg U$ are equivalent to the statement that $U = O(V)$
for positive functions $U$ and $V$). We also use the symbol `$o$'
with its usual meaning: the statement $U=o(V)$ is equivalent to
$U/V\to 0$.

We always use the letters $p$ and $q$ to denote prime numbers.

\section{Divisors of the Ramanujan function}

\begin{theor}[\!]
There exist infinitely many $n$ such that \hbox{$\tau(n) \neq 0$}
and $P(\tau(n)) \geq (\log n)^{33/31 + o(1)}$.
\end{theor}

\begin{proof}
For a constant  $A>0$ and a real $z$ we
define the set
\begin{equation*}
\cS_A(z)=\{n\leq z\hbox{:}\ P(n)\leq (\log n)^A\}.
\end{equation*}
For every $A>1$, we have $\#\cS_A(z)= z^{1-1/A+o(1)}$,
as $z \to \infty$ (see eq.~(1.14) in~\cite{HiTe} or Theorem~2 in \S~III.5.1
of~\cite{Te}).

Let $x> 0$ be sufficiently large. By a result of
Serre~\cite{Serre}, the estimate $\#\{p\leq y\hbox{:}\
\tau(p)=0\}\ll y/(\log y)^{3/2}$ holds as $y$ tends to infinity.
Applying this estimate with $y=x^{1/2}$, it follows that there are
only $ o( \pi(y))$ primes $p<y$ such that $\tau(p) = 0$. It is
also obvious from~\eqref{eq:quadr}  that $\tau(p^2) \neq 0$.

Assume that for some $A$ with $1 < A < 33/31$, we have the
inequality $P(\tau(p)\tau(p^2)) \leq (\log y)^A$ for all remaining
primes \hbox{$p \leq y$}. We see from~\eqref{eq:quadr}
and~\eqref{eq:Deligne} that $|\tau(p^2)|=|\tau(p)^2-p^{11}|\leq
3p^{11} \leq 3y^{11}$. Denoting \hbox{$z_1 =3y^{11}$} and
\hbox{$z_2 = 2 p^{11/2}$}, we deduce  that for  $( 1+ o(1)) \pi(y)
= y^{1 + o(1)}$ primes $p< y$ with $\tau(p) \neq 0$, we have a
representation $p^{11} =  s_1^2 - s_2$, where $s_i \in
\cS_A(z_i)$, $i=1,2$. Thus
\begin{equation*}
y^{1 + o(1)} \leq \# \cS_A(z_1) \# \cS_A(z_2)  \leq (
z_1z_2)^{1-1/A+o(1)} \leq (6 y^{33/2})^{1-1/A+o(1)},
\end{equation*}
which is impossible for $A < 33/31$. This  completes the proof.\hfill $\Box$
\end{proof}

We remark in passing that the above proof shows that the inequality
$P(\tau(p)\tau(p^2))>(\log p)^{33/31+o(1)}$ holds for almost all primes
$p$.

\begin{theor}[\!]
\label{thm:omega}
The estimate
\begin{equation*}
\omega\left(\prod_{\substack{p<x^{1/3}\\ \tau(p)\neq 0}}\tau(p)\tau(p^2)
\tau(p^3)\right)\geq \left(\frac{1}{6\log 7}+o(1)\right) \log x
\end{equation*}
holds as $x$ tends to infinity.
\end{theor}

\begin{proof}
Let $x$ be a large positive integer and put $y=x^{1/3}$.
Let $\cR$ be   the  set of   odd primes \hbox{$p\leq y$} such
that $\tau(p)\neq 0$. Note that since $\tau(p)\neq 0$, it follows
that $\tau(p^2)\neq 0$ and $\tau(p^3)\neq 0$. Let
\begin{equation*}
M=\prod_{p\in \cR}\tau(p)\tau(p^2)\tau(p^3)\mand
s=\omega(M).
\end{equation*}
Since $\tau(p^2)=\tau(p)^2-p^{11}$ and $\tau(p^3)=\tau(p)(\tau(p)^2-2p^{11})$,
eliminating $p^{11}$, we get the equation
\begin{equation*}
1=\frac{2\tau(p^2)}{\tau(p)^2}-\frac{\tau(p^3)}{\tau(p)^3}.
\end{equation*}
We claim that the rational numbers $2\tau(p^2)/\tau(p)^2$ are
distinct for distinct odd primes. Indeed, if
$\tau(p_1^2)/\tau(p_1)^2= \tau(p_2^2)/\tau(p_2)^2$ for two
distinct odd  primes $p_1, p_2$, we get that
$p_1^{11}/\tau(p_1)^2=p_2^{11}/ \tau(p_2)^2$, or  $p_1^{11}
\tau(p_2)^2=p_2^{11} \tau(p_1)^2$. Therefore,
$p_1^{11}|\tau(p_1)^2$. Thus, $p_1^{12}|\tau(p_1)^2$, which is
impossible for $p_1>3$ because of~\eqref{eq:Deligne}, and can be
checked by hand to be impossible for $p_1=3$.

Let $\cS$ be the set of all prime divisors of $M$. Thus, $\# \cS =
s$. We see that the equation $u-v=1$ has $\# \cR$ distinct
solutions in the $\cS$-units
\begin{equation}\label{eq:Solution}
(u,v)=\left(\frac{2\tau(p^2)}{\tau(p)^2},
\frac{\tau(p^3)}{\tau(p)^3}\right).
\end{equation}
It is known (see~\cite{Eve}), that the number of solutions of such
a $\cS$-unit equation is $O(7^{2s})$. We thus get that $7^{2s}\gg
\#\cR =  (1+o(1))\pi(y)$, giving
\begin{equation*}
s\geq \frac{1}{6\log 7}(1+o(1))\log x
\end{equation*}
as $x\rightarrow \infty$, which finishes the proof.\hfill $\Box$
\end{proof}

\begin{theor}[\!]
\label{thm:TripleProd}
The estimate
\begin{equation*}
P(\tau(p)\tau(p^2)\tau(p^3))>(1+o(1))\frac{\log\log p\log\log\log p}
{\log\log\log\log p}
\end{equation*}
holds as $p$ tends to infinity through primes such that
$\tau(p)\neq 0$.
\end{theor}

\begin{proof}
As in the proof of Theorem~\ref{thm:omega}, we consider the equation
$u-v=1$,
having  the solution~\eqref{eq:Solution} for every prime $p$ with
$\tau(p) \neq 0$. Write
\begin{equation*}
u = E/D \mand v = F/D,
\end{equation*}
where $D$  is the smallest positive common denominator of $u$ and $v$.
Then
\begin{equation*}
E = Du=2D-2p^{11}D/\tau(p)^2 \mand F = Dv=D-2Dp^{11}/\tau(p)^2
\end{equation*}
are integers with $\gcd(E,F) = 1$, and since $E-F = D$, we
also have  $\gcd(D,E) = \gcd(D,F) = 1$.

We note the inequalities
\begin{equation}
\label{eq:Inequal}
D \ll p^{11} \mand p \ll \max\{|E|,|F|\}\ll p^{22} .
\end{equation}
Indeed, the upper bounds follow directly from~\eqref{eq:Deligne}.
It also follows from~\eqref{eq:Deligne} that
\hbox{$p^6\hbox{{$\not|$}}~\tau(p)$}. This shows that
$p^{11}/\tau(p)^2$ is a rational number whose numerator is a
multiple of $p$. In particular,
\begin{equation*}
E - 2 F  = \frac{2Dp^{11}}{\tau(p)^2} \geq p,
\end{equation*}
which implies the lower bound in~\eqref{eq:Inequal}.\pagebreak

We have \hbox{$P(\tau(p)\tau(p^2)\tau(p^3)) \geq \ell$}, where \hbox{$\ell = P(EDF)$}.

Let $t = \omega(\tau(p)\tau(p^2)\tau(p^3))$. By~\eqref{eq:Inequal}, we
see that there exists a prime $q$ and a positive integer $\alpha$ such
that $q^{\alpha}$ divides one of $E$ or $F$ and $q^{\alpha}\gg p^{1/t}$.

First we assume that $q^{\alpha}|E = D-F$, and  write
\begin{equation*}
D=\prod_{j=1}^t q_j^{\beta_j} \mand
F=\prod_{j=1}^t q_j^{\gamma_j},
\end{equation*}
with some primes $q_j$ and non-negative integers $\beta_j,\gamma_j$ such
that $\min\{\beta_j,\gamma_j\}=0$ for all $j=1,\dots, t$ (clearly,
$\beta_i =\gamma_i=0$ for $q_i = q$). By~\eqref{eq:Inequal}, we also
have
\begin{equation*}
B = \max_{j=1, \dots, t} \{\beta_j,\gamma_j\}\ll \max\{\log D, \log
|E|\} \ll \log p.
\end{equation*}
Using
the lower bound for linear forms in $q$-adic logarithms  of Yu~\cite{Yu},
we derive
\begin{equation}
\label{eq:MainIneq}
\alpha \leq q  c^t \log B \prod_{j=1}^t \log q_j \ll
    \ell   ( c \log \ell )^{t} \log\log p
\end{equation}
with some absolute constant $c > 0$.
Since also
\begin{equation*}
\alpha \gg \frac{\log p}{t \log q} \geq  \frac{\log p}{t\log \ell },
\end{equation*}
we get
\begin{equation*}
\frac{\log p}{\log\log p}\ll \ell t  ( c \log \ell )^{t}  \ll
\ell   ( 2 c \log \ell )^{t}.
\end{equation*}
Hence,
\begin{equation}
\label{eq:tl}
\log\log p\leq t(1+o(1))\log\log\ell.
\end{equation}
By the prime number theorem (see \cite{HardyWright}), we have
\begin{equation*}
t \leq (1  + o(1))  \frac{\ell }{\log \ell},
\end{equation*}
which together with \eqref{eq:tl} leads us to
\begin{equation*}
(1+o(1))\frac{\log\log p\log\log\log p}
{\log\log\log\log p}\leq t.
\end{equation*}

The case $q^{\alpha}|F = D-E$ can be considered completely analogously
which concludes the proof.\hfill $\Box$
\end{proof}

We recall that the $ABC$-conjecture asserts that
for any fixed $\varepsilon > 0$ the inequality
\begin{equation*}
Q(abc) \gg (\max{|a|, |b|, |c|})^{1- \varepsilon}
\end{equation*}
holds for any relatively prime integers $a,b,c$ with $a + b = c$.
Thus, in the notation of the proof of Theorem~\ref{thm:TripleProd},
we immediately conclude from~\eqref{eq:Inequal} that the
$ABC$-conjecture yields
\begin{equation*}
Q(\tau(p)\tau(p^2)\tau(p^3))\geq Q(DEF) \geq  p^{1 + o(1)}.
\end{equation*}
Thus, by the prime number theorem,
\begin{equation*}
P(\tau(p)\tau(p^2)\tau(p^3)) \geq (1 + o(1)) \log p.
\end{equation*}
The best known unconditional result of Stewart  and Yu~\cite{SteYu}
towards the $ABC$-conjecture implies that
\begin{equation*}
Q(\tau(p)\tau(p^2)\tau(p^3)) \geq Q(DEF) \geq  (\log p)^{3 + o(1)}.
\end{equation*}

\section{Factorials and the Ramanujan function}

In~\cite{Lu}, all the positive integer solutions $(m,n)$ of the equation
$f(m!)=n!$ were found, where $f$ is any one of the multiplicative
arithmetical functions $\varphi$, $\sigma$, $d$, which are the Euler
function, the sum of divisors function, and the number of divisors
function, respectively. Further results on such problems have been
obtained by Baczkowski~\cite{Ba}. Here, we study this problem for the
Ramanujan function.

\begin{theor}[\!]
There are only finitely many effectively computable pairs of positive integers
$(m,n)$ such that $|\tau(m!)|=n!$.
\end{theor}

\begin{proof}
Assume that $(m,n)$ are positive integers such that $\tau(m!)=n!$.
By~\eqref{eq:Deligne} and the Stirling formula
\begin{align*}
\exp((1+o(1)) n \log n) & = n!=\tau(m!)<(m!)^{11/2+o(1)}\\[.3pc]
& < \exp((11/2+o(1)) m \log m),
\end{align*}
as $m$ tends to infinity. Thus, we conclude that if $m$ is sufficiently
large, then $n<6m$.

Let $\nu(m)$ be the order at which the prime $2$ appears in the
prime factorization of $m!$. It is clear that $\nu(m)>m/2$ if $m$
is sufficiently large. Since $\tau$ is multiplicative, it follows
that $u_{\nu(m)}=\tau(2^{\nu(m)})|n!$, where the Lucas sequence
$u_r$ is given by~\eqref{eq:Bin Rec} with $u_0 = 1$, $u_1 = -24$.

For \hbox{$r\geq 1$}, we put $\zeta_r=\exp(2\pi i/r)$ and consider the
sequence $v_r=\Phi_r(\alpha,\beta)$ where
\begin{equation*}
\Phi_r(X,Y) = \prod_{\substack{1\leq k\leq r\\ \gcd(k,r)=1}}(X-\zeta_r^k Y).
\end{equation*}
It is known that $v_r | u_r$. It is also known (see~\cite{BHV}),
that $v_r=A_rB_r$, where $A_r$ and $B_r>0$ are integers,
\hbox{$|A_r|\leq 6(r+1)$} and every prime factor of $B_r$ is
congruent to \hbox{$\pm 1\!\!\!\pmod{r+1}$}. Let $\alpha$ and $\beta$
be the two roots of the characteristic equation \hbox{$\lambda^2 -
24\lambda - 2048=0$}. Since both inequalities \hbox{$|v_k|\leq
2|\alpha|^{k+1}$} and \hbox{$|v_k|\geq |\alpha|^{k+1-\gamma\log
(k+1)}$} hold for all positive integers $k$ with some absolute
constant $\gamma$ (see, for example, Theorem~3.1 on p.~64 in
\cite{ST}), it follows that
\begin{align*}
6(r+1)B_r &\geq 2^{-\tau(r+1)}\alpha^{\varphi(r+1)-\gamma \tau(r+1)\log (r+1)}\\[.3pc]
&= |\alpha|^{\varphi(r+1)+O(\tau(r+1)\log (r+1))}.
\end{align*}
Since $\varphi(r+1)\gg r/\log\log r$, and
$\tau(r+1)\log(r+1)=r^{o(1)}$, the above inequality implies that
\begin{equation*}
B_r>|\alpha|^{\varphi(r+1)/2}
\end{equation*}
whenever $r$ is sufficiently large.\pagebreak

In particular, we see that $B_{\nu(m)} | \tau(m!)$, has all prime
factors \hbox{$\ell \equiv \pm 1\!\!\!\pmod{\nu(m)+1}$}, and is of the size
\begin{equation*}
B_{\nu(m)} >\exp(c m/\log\log m),
\end{equation*}
where $c$  is some positive constant.

However, since $B_{\nu(m)}|n!$ and $n<6m$, it follows that all
prime factors $\ell$ of $B_{\nu(m)}$ satisfy $\ell < 6m$. Since
$\nu(m)>m/2$, there are at most $26$ primes $\ell <6m$ with $\ell
\equiv \pm 1\!\pmod{\nu(m)+1}$. Furthermore, again since
$B_{\nu(m)}|n!$, $n<6m$, and all prime factors $\ell$ of
$B_{\nu(m)}$ satisfy \hbox{$\ell \equiv \pm 1\!\!\!\pmod{\nu(m)+1}$}, it
follows that $\ell^{14} \nmid B_{\nu(m)}$. Hence,
\begin{equation*}
B_{\nu(m)}<(6m)^{26\cdot 13}= m^{O(1)}.
\end{equation*}
Comparing this with the above lower bound on $B_{\nu(m)}$, we conclude
that  $m$ is bounded.

\hfill $\Box$
\end{proof}

\section*{Acknowledgements}

During the preparation of this paper, the first author was
supported in part by grants SEP-CONACYT 37259-E and 37260-E, and
the second author was supported in part by ARC grant DP0211459.


\begin{thebibliography}{99}
\bibitem{Ba} Baczkowski~D, Master Thesis (Miami Univ., Ohio, 2004)

\bibitem{BHV} Bilu~Y, Hanrot~G and Voutier~P~M, Existence of
primitive divisors of Lucas and Lehmer numbers, with an appendix
by M~Mignotte, {\it J.~Reine Angew. Math.} {\bf 539} (2001) 75--122

\bibitem{Eve} Evertse~J-H, On equations in $S$-units and the
Thue-Mahler equation, {\it Invent. Math.} {\bf 75} (1984) 561--584

\bibitem{HardyWright} Hardy~G~H and Wright~E~M, {An introduction
to the theory of numbers} (Oxford Univ. Press, Oxford, 1979)

\bibitem{HiTe} Hildebrand~A and Tenenbaum~G, Integers without large
prime factors, {\it J. de Th{\'e}orie des Nombres de Bordeaux} {\bf
5} (1993) 411--484

\bibitem{Lu} Luca~F, Equations involving arithmetic functions of
factorials, {\it Divulg. Math.} {\bf 8(1)} (2000) 15--23

\bibitem{Murty} Murty~M~R, The Ramanujan $\tau$ function,
Ramanujan revisited, {\it Proc. Illinois Conference on Ramanujan} (1988)
269--288

\bibitem{MM} Murty~M~R and Murty~V~K, Prime divisors of Fourier
coefficients of modular forms, {\it Duke Math. J.} {\bf 51} (1985)
521--533

\bibitem{MMS} Murty~M~R, Murty~V~K and Shorey~T~N, Odd values of
the Ramanujan $\tau$-function, {\it Bull. Soc. Math. France} {\bf
115} (1987) 391--395

\bibitem{MuWo} Murty~M~R and Wong~S, The $ABC$ conjecture and prime
divisors of the Lucas and Lehmer sequences, Number Theory for the
Millennium, vol.~III (MA: A.~K.~Peters, Natick) (2002) 43--54

\bibitem{Serre} Serre~J~P, Quelques applications du th{\'e}or{\`e}me
de densit{\'e} de Chebotarev, {\it Publ. Math., Inst. Hautes {\'E}tud.
Sci.} {\bf 54} (1981) 123--201

\bibitem{ST} Shorey~T~N and Tijdeman~R, Exponential diophantine
equations (Cambridge: Cambridge Univ. Press) (1986)

\bibitem{Shorey} Shorey~T~N, Ramanujan and binary recursive
sequences, {\it J.~Indian Math. Soc.} {\bf 52} (1987) 147--157

\bibitem{Stew} Stewart~C~L, On divisors of Fermat, Fibonacci, Lucas and
Lehmer numbers, III, {\it J.~London Math. Soc.} {\bf 28} (1983)
211--217

\bibitem{SteYu} Stewart~C~L and Yu~K~R, On the $abc$ conjecture,
II, {\it Duke Math.~J.} {\bf 108} (2001) 169--181

\bibitem{Te} Tenenbaum~G, Introduction to analytic and
probabilistic number theory (Cambridge: Cambridge Univ. Press) (1995)

\bibitem{Yu} Yu~K, $p$-Adic logarithmic forms and group varieties,
II, {\it Acta Arith.} {\bf 89} (1999) 337--378
\end{thebibliography}
\end{document}